\documentclass[10pt,twoside]{article}
\usepackage{amsmath,amsbsy,amsfonts,amsthm,amssymb}
\usepackage{color}
\newtheorem{lemma}{Lemma}[section]
\newtheorem{proposition}{Proposition}[section]
\newtheorem{theorem}{Theorem}[section]
\newtheorem{corollary}{Corollary}[section]
\newtheorem{remark}{Remark}[section]

\let\Section=\section
\def\section{\setcounter{equation}{0}\Section}

\title{Singularly perturbed biharmonic problems with  superlinear nonlinearities}

\author{Marcos T. O. Pimenta\thanks{Research Supported by CNPq - Brazil}\\
\noindent Departamento de Matem\'{a}tica\\
\noindent Centro de Ci\^{e}ncias Exatas\\
\noindent Universidade Estadual de Londrina - UEL \\
\noindent 86051-980, Londrina - PR, Brazil.\\
\noindent e-mail: {\tt{mtopimenta@uel.br}}\\
\mbox{}\\
and \\
\mbox{}\\
\noindent S\'{e}rgio H. M. Soares\thanks{S\'{e}rgio H. M. Soares was partially supported by CNPq - Brazil} \\
\noindent Departamento de Matem\'atica \\
\noindent Instituto de Ci\^{e}ncias Matem\'{a}ticas e de Computa\c{c}\~{a}o\\
\noindent Universidade de S\~ao Paulo \\
\noindent 13560-970, S\~ao Carlos - SP, Brazil. \\
\noindent e-mail: {\tt{monari@icmc.usp.br}}\\
}
\date{}

\begin{document}

\maketitle

\begin{abstract}
We are interested in finding  a family of solutions to a singularly perturbed  biharmonic equation which has a concentration behavior. The proof is based on  variational methods  and it is used a weak version of the Ambrosetti-Rabinowitz condition.
\end{abstract}


\section{Introduction}

This paper was motivated by some results for the following class of semilinear elliptic equations
\begin{equation}
\left\{
\begin{array}{l}
\epsilon^2\Delta u + V(x)u = f(u), \ \ \mbox{in
$\mathbb{R}^N$},\\
u\in H^1(\mathbb{R}^N).
\end{array} \right. \label{P1}
\end{equation}
This problem has recently been extensively studied, see for example \cite{Alves2,  Del Pino, Floer, Jeanjean, Oh1, Oh2, Rabinowitz, Wang} and the references therein.   The existence and concentration of spike-layered solutions was first studied by Floer and Weinstein in \cite{Floer} in the one dimensional case.  Later, Oh in \cite{Oh1} and \cite{Oh2} extended this result to higher dimensions considering a larger class of nonlinearities. These results have inspired Rabinowitz in \cite{Rabinowitz} to deal with this class of problems, considering the so called Rabinowitz condition under the potential V,
$$0 < V_0:=\inf_{\mathbb{R}^N}V < \liminf_{|x| \to \infty} V(x) =: V_\infty.$$
In \cite{Rabinowitz}, it is used a  mountain-pass type argument to show the existence of a ground-state solution to (\ref{P1}) where $\epsilon = 1$.    In \cite{Wang}, Wang proves that the maximum points of the solutions obtained in \cite{Rabinowitz}  converge to a global minimum point of $V$ as $\epsilon \rightarrow 0$, characterizing the concentration behavior of this family of solutions.  In  \cite{Del Pino}, del Pino and Felmer developed a method to obtain a family of solutions concentrating around a local minimum point of $V$.   In \cite{Jeanjean}, Jeanjean and Tanaka proved the same result obtained in \cite{Del Pino}, but with the nonlinearity $f$ satisfying weaker assumptions.  More specifically, they considered the case where $f$  neither satisfies   the monotonicity condition on the function $s \mapsto f(s)/s$, nor  the so called Ambrosetti-Rabinowitz condition
\begin{itemize}
\item [$(AR)$] $0 < \mu F(s) \leq f(s)s,$ \quad for all $s \neq 0$ and for some $\mu > 2$.
\end{itemize}
The purpose of this paper is to provide similar results to the following biharmonic Schr\"{o}dinger elliptic equation
\begin{equation}
\left\{
\begin{array}{l}
\epsilon^4\Delta^2u + V(x)u = f(u), \ \ \mbox{in
$\mathbb{R}^N$},\\
u\in H^2(\mathbb{R}^N).
\end{array} \right. \label{P2}
\end{equation}
The nonlinearity $f$ will be assumed to satisfy a weaker superlinearity condition than $(AR)$. More specifically, we assume the following  conditions on $f$ and $V$:
\begin{description}
\item [$(V_1)$] $V\in C^0(\mathbb{R}^N) \cap L^{\infty}(\mathbb{R}^N)$.
\item [$(V_2)$] $0 < V_0:=\displaystyle\inf_{\mathbb{R}^N}V < \displaystyle\liminf_{|x| \to \infty} V(x) =: V_\infty$.
\item [$(f_1)$] $f\in C^1(\mathbb{R})$.
\item [$(f_2)$] $f(0)=f'(0)=0$.
\item [$(f_3)$] There exist $c_1,c_2>0$ and
$p\in\left(1,2_* - 1\right)$ such that $|f(s)|\leq c_1 |s| +
c_2|s|^p$ for all $s\in\mathbb{R}$, where $2_* = {2N}/{(N-4)}$.
\item [$(f_4)$] $\displaystyle\lim_{|s|\rightarrow \infty}\frac{F(s)}{s^2} = +\infty$, where $F(s) = \int_0^s f(t)dt$.
\item [$(f_5)$] $\displaystyle\frac{f(s)}{s}$ is increasing for $s > 0$ and decreasing for $s < 0$.
\end{description}

\begin{remark} \label{remark1}
The conditions $(f_2)$ ad $(f_5)$ imply that
$$
H(s):=f(s)s - 2F(s) >0,\quad F(s) >0, \mbox{ and}\quad sf(s) >0,\quad   \forall\, s \neq 0.
$$
\end{remark}
Our main result is the following.

\begin{theorem}
Assume that conditions \emph{$(V_1)$}, \emph{$(V_2)$} and \emph{$(f_1) - (f_5)$} hold. Then for each sequence $\epsilon_n\rightarrow 0$, along a subsequence, there exists a nontrivial weak solution $u_n$ of (\ref{P2}) (with $\epsilon=\epsilon_n$). Moreover, if $x_n$ is the maximum point of $|u_n|$, then $$\lim_{n\to\infty}V(x_n) = \inf_{\mathbb{R}^N} V.$$ \label{theorem1.1}
\end{theorem}
In  \cite{Pimenta}, we establish the same conclusion of Theorem  \ref{theorem1.1} in the case of the potential $V$ satisfies a local condition given by del Pino and Felmer in  \cite{Del Pino}.

 Although so many of our arguments were inspired in the works mentioned above, it is worth pointing out that some of them have to be deeply modified because of some difficulties that the lack a general maximum principle to the biharmonic operator gives rise. For instance, in \cite{Wang} Wang uses a Harnack type inequality to prove the uniform decay of some translations of solutions that we were not able to find to biharmonic subsolutions. Hence, we use an $L^\infty$ estimate from Ramos \cite{Ramos} and an $L^p$ estimate from Agmon \cite{Agmon} in order to prove the same result to the fourth-order operator. Some arguments about compactness in Nehari manifolds found in \cite{Alves1} seems to be useful in this argument too. Finally, the lack of a standard form of the Ambrosetti-Rabinowitz condition in our work represents some difficulty to prove that the $(PS)$ sequences are bounded, which required some arguments of Miyagaki and Souto in \cite{Souto} and also represents a difficulty to prove that the Nehari manifold is homeomorphic to the unitary sphere in $H^2(\mathbb{R}^N)$. This last problem can be dropped out using some arguments of Weth and Szulkin in \cite{Weth}.

%it is here employed similar methods to those above-mentioned,  many arguments have to be deeply modified because of the lack of a general maximum principle to the biharmonic %operator. Another difficulty in dealing with biharmonic equations is the application of  Moser's iteration technique which would be useful in our method. In order to overcome these d%ifficulties, we use some compactness results on Nehari manifolds, some arguments found in \cite{Alves6} and an a priori estimate for solutions of subcritical biharmonic problems found i%n \cite{Ramos}.

This paper is organized as follows. In the second section, we use some arguments of \cite{Rabinowitz} to prove the existence of a family of solutions to (\ref{P2}).  The third section is devoted to prove that this family has a concentration behavior.

\section{Existence}

\noindent In main result this section is the following:

\begin{theorem}
Let assumptions $(V_1),(V_2)$ and $(f_1) - (f_5)$ hold. Then there exists $\epsilon_0 > 0$ such that problem (\ref{P2}) has a nontrivial weak solution $u_\epsilon$ provided that $\epsilon < \epsilon_0$.
\label{theorem2.1}
\end{theorem}
We observe that (\ref{P2}) is equivalent  to the problem
\begin{equation}
\left\{
\begin{array}{l}
\Delta^2v + V(\epsilon x)v = f(v), \ \ \mbox{in
$\mathbb{R}^N$},\\
u\in H^2(\mathbb{R}^N),
\end{array} \right. \label{P3}
\end{equation}
and the  equivalence among the solutions $u_\epsilon$ of (\ref{P2}) and $v_\epsilon$ of (\ref{P3}) is given by $u_\epsilon(\epsilon x) = v_\epsilon(x)$.

In order to use variational methods, lwe consider the Sobolev space $H^2(\mathbb{R}^N)$ endowed with the inner product
$$\langle u ,v \rangle_\epsilon = \int_{\mathbb{R}^N}\left(\Delta u \Delta v + V(\epsilon x)uv \right) dx,$$
which gives rise to the following norm
$$\|u\|_\epsilon = \left(\int_{\mathbb{R}^N}\left(|\Delta u|^2 + V(\epsilon x)u^2 \right) dx\right)^{1/2}.$$
From now on we denote by $E_\epsilon = \left( H^2(\mathbb{R}^N),\langle \cdot , \cdot \rangle_\epsilon\right)$.  We consider the functional  $I_\epsilon$ defined on  $E_\epsilon$ by
$$
I_\epsilon(u) = \frac{1}{2}\int_{\mathbb{R}^N}\left(|\Delta u|^2 + V(\epsilon x)u^2\right)dx - \int_{\mathbb{R}^N}F(u)dx,$$
where $F(s)=\int_0^s f(t)dt$.  The functional $I_\epsilon \in C^1(E_\epsilon, \mathbb{R}^N)$ and
\[
I'_\epsilon(u)v = \int_{\mathbb{R}^N}\left(\Delta u \Delta v + V(\epsilon x)uv \right) dx - \int_{\mathbb{R}^N}f(u)vdx,
\]
for all $u, v \in E_\epsilon$. Hence, critical points of  $I_\epsilon$ are weak solutions of (\ref{P3}).

Our first lemma provides conditions under which $I_\epsilon$ satisfies the geometric hypotheses of the  Mountain Pass Theorem.

\begin{lemma}\label{lemma1}
Assume that conditions $(f_2) - (f_4)$ hold. Then, for each $\epsilon > 0$ there exist $\rho, r>0$ and $\varphi\in E_\epsilon$ with $\|\varphi\|_\epsilon > r$, such that
\begin{itemize}
\item[i)] $I_\epsilon(u) \geq \rho$ for all $\|u\|_\epsilon = r$;
\item[ii)] $I_\epsilon(\varphi) < 0$.
\end{itemize}
\label{lemma2.1}
\end{lemma}
\noindent \textbf{Proof.}
Using $(f_2)$ and $(f_3)$ and the Sobolev embeeding, we can prove that for all $\eta > 0$, there exists a constant $C(\eta) > 0$ such that
\[
\int_{\mathbb{R}^N}|F(u)|dx \leq
\eta\|u\|_\epsilon^2 + C(\eta)\|u\|_\epsilon^{p+1}.
\]
Hence, by choosing $\eta \in (0, 1/2)$, there exists a small $r>0$ in such a way that
\[
I_\epsilon (u) \geq \rho > 0,\quad \mbox{for all $\|u\|_\epsilon=r$},
\]
where $\rho = [(1/2 - \eta) - C(\eta)r^{p-1}]r^2$. This establishes {\it i)}.

In order to prove {\it ii)}, fix $\varphi\in C^\infty_0(\mathbb{R}^N)$ with $\varphi > 0$.  By $(f_4)$,  for every  $M \geq {\|\varphi\|_\epsilon^2}/{2\|\varphi\|_{L^2}^2}$, there exists a constant $c_0 > 0$ such that
$$F(s) \geq M|s|^2 - c_0, \quad \mbox{for all $s\in \mathbb{R}$}.$$
Then,
\begin{eqnarray*}
I_\epsilon(t\varphi) & = & \frac{t^2}{2}\|\varphi\|_\epsilon^2 - \int_{\mathbb{R}^N}F(t\varphi)dx\\
& \leq & \frac{t^2}{2}\|\varphi\|_\epsilon^2 - t^2 M\int_{\mathbb{R}^N}|\varphi|^2
dx + c_0|supp(\varphi)|\\
& = & t^2 \left(\frac{\|\varphi\|_\epsilon^2}{2} - M\int_{\mathbb{R}^N}|\varphi|^2dx\right) + c_0|supp(\varphi)|.
\end{eqnarray*}
Therefore, $I_\epsilon(tu)\rightarrow -\infty$ as $t\rightarrow +\infty$ and the proof is
complete. \quad\hfill $\Box$

As a consequence of Lemma \ref{lemma2.1}, for every $\epsilon > 0$,there corresponds a minimax  value associated with  (\ref{P3}) and given by
\begin{equation}
c_\epsilon = \inf_{g\in\Gamma_\epsilon}\sup_{0 \leq t \leq
1}I_\epsilon(g(t)), \label{c_epsilon}
\end{equation}
where $$\Gamma_\epsilon = \left\{ g\in C([0,1],E_\epsilon); \,\
g(0) = 0 \,\ \mbox{e} \,\ I_\epsilon(g(1))<0\right\}.$$
In order to get least energy solutions for (\ref{P3}), consider the Nehari manifold
\begin{equation}
\mathcal{N}_\epsilon=\left\{u\in E_\epsilon\backslash\{0\};\ \
I_\epsilon'(u)u=0\right\}. \label{nehari}
\end{equation}
Unlike in \cite{Rabinowitz}, when $f$ does not satisfies the  Ambrosetti-Rabinowitz condition things become much more difficult to prove that $\mathcal{N}_\epsilon$ is homeomorphic to the unitary sphere in $E_\epsilon$. However, following some arguments found in \cite{Weth}, we can show  that $\mathcal{N}_\epsilon$ is homeomorphic to the unitary sphere using the superlinearity condition $(f_4)$. Hence, similar analysis to that in  \cite{Rabinowitz} shows that
$$
c_\epsilon = \inf_{u\in E_\epsilon\backslash\{0\}}\max_{t\geq
0}I_\epsilon(tu) = \inf_{u \in \mathcal{N}_\epsilon}I_\epsilon(u)$$
and
$\mathcal{N}_\epsilon = \{\varphi_\epsilon(u)u;\: u\in E_\epsilon\backslash\{0\} \}$, where $\varphi_\epsilon(u) > 0$ is such that $I_\epsilon(\varphi_\epsilon(u) u) = \max_{t \geq 0}I_\epsilon(tu)$. Hence,  every solution in the level $c_\epsilon$  is a least energy solution.

We now use some arguments of Jeanjean and Tanaka in \cite{Jeanjean} and Miyagaki and Souto in \cite{Souto} to prove that the Palais-Smale sequences of $I_\epsilon$ are bounded.

\begin{lemma}
If $(v_n)$ is a $(PS)_c$ sequence for $I_\epsilon$, then $(v_n)$ is bounded in $E_\epsilon$.
\label{lemma2.2}
\end{lemma}
\noindent \textbf{Proof.}  Suppose by contradiction that $\|v_n\|_\epsilon \rightarrow \infty$ as $n \rightarrow \infty$. Let us define $$w_n = \frac{v_n}{\|v_n\|_\epsilon}.$$
We claim that one of the two statements holds:
    \begin{description}
    \item [{\it i)}] $w_n \rightarrow 0$ in $L^r(\mathbb{R}^N)$, for all $2  <  r  <  2_*$.
    \item [{\it ii)}] There exist $(y_n)\subset \mathbb{R}^N$ and constants $R,\beta>0$ such that
    $$\liminf_{ n\to \infty} \int_{B_R(y_n)}w_n^2dx \geq \beta.$$
    \end{description}
    Indeed, suppose that {\it ii)} does not hold. Then for all $R>0$,
    $$\sup_{y \in \mathbb{R}^N} \int_{B_R(y)}w_n^2 dx \rightarrow 0, \quad \mbox{as $n\to\infty$}.$$ Once $(w_n)$ is bounded in $L^2(\mathbb{R}^N)$ and $(\nabla w_n)$ is bounded in $L^{2^*}(\mathbb{R}^N)$, where $2^* = 2N/(N-2)$,  Lions Lemma  \cite{Lions}  implies that $w_n \rightarrow  0$ in $ L^r(\mathbb{R}^N)$, for all $2 < r < N2^*/(N-2^*) = 2_*$. Hence {\it (i)} is proved to hold, and the claim is verified.

 Suppose that  {\it i)} holds. By $(f_2)$, $(f_3)$ and {\it i)}, for all $\mu>0$,
    \begin{equation}
    \lim_{n\to\infty} \int_{\mathbb{R}^N}F(\mu w_n)dx = 0.
    \label{eq3}
    \end{equation}
   Set $s_n\in [0,1]$ such that
    $$
I_\epsilon(s_nv_n) = \max_{t\in [0,1]}I_\epsilon(tv_n).
$$
For every $\mu >0$ and  $n$ sufficiently large, we have
    \begin{equation}
    I_\epsilon(s_nv_n) \geq I_\epsilon\left(\frac{\mu}{\|v_n\|_\epsilon}v_n\right) = \frac{\mu^2}{2} - \int_{\mathbb{R}^N}F(\mu w_n)dx.
    \label{eq4}
    \end{equation}
    By (\ref{eq3}) and (\ref{eq4}) it follows that
    $$\liminf_{n\to\infty} I_\epsilon(s_nv_n) \geq \frac{\mu^2}{2}, \quad \mbox{for all $\mu>0$.}$$ Hence
    \begin{equation}
    \liminf_{n\to\infty} I_\epsilon(s_nv_n) = +\infty.
    \label{eq5}
    \end{equation}
Since $I_\epsilon(0) = 0$ and $I_\epsilon(v_n) \rightarrow c$ as $n\to\infty$, we have $s_n \in (0,1)$ for every $n$ sufficiently large. Therefore,  $I'_\epsilon(s_nv_n)s_nv_n = 0$.  By Remark \ref{remark1}, for all $t \in [0,1]$, it follows that
    \begin{eqnarray*}
    2I_\epsilon(tv_n)  \leq  2I_\epsilon(s_nv_n) - I_\epsilon'(s_nv_n)s_nv_n
     \leq   2I_\epsilon(v_n) + o_n(1) \leq C_1 + o_n(1).
    \end{eqnarray*}
Given any  $R_0 > 0$, there exists $n_0> 0$ such that  ${R_0}/{\|v_n\|_\epsilon} < 1$ for all $n\geq n_0$. Hence,
    \begin{equation}
    2I(R_0w_n) = 2I_\epsilon\left(\frac{R_0}{\|v_n\|}v_n\right) \leq C_1 + o_n(1). \label{eq6}
    \end{equation}
 On the other hand,
$$2I_\epsilon(R_0w_n) = R_0^2 - 2\int_{\mathbb{R}^N} F(R_0w_n)dx  = R_0^2 + o_n(1),$$
which contradicts (\ref{eq6}) because $R_0 > 0$ is arbitrary.  Therefore, we conclude that {\it i)} does not hold, and consequently {\it ii)} occurs. From this, we can define $\bar{w}_n(x) := w_n(x + y_n)$. As $(w_n)$, ($\bar{w}_n)$ is bounded in $H^2(\mathbb{R}^N)$ as well. Hence, there exists $\bar{w} \in H^2(\mathbb{R}^N)$ such that $\bar{w}_n \rightharpoonup \bar{w}$ in $H^2(\mathbb{R}^N)$ along a subsequence. Moreover, by {\it ii)}, $\bar{w} \neq 0$ in $H^2(\mathbb{R}^N)$  and $\bar{w}(x) \neq 0$ almost everywhere in a subset $\Sigma$ of $B_R(0)$ with positive measure. Since $(I_\epsilon(v_n))$ is bounded, we have
    $$\frac{1}{2} + o_n(1) = \int_{\mathbb{R}^N}\frac{F(v_n)}{\|v_n\|^2_\epsilon}dx = \int_{\mathbb{R}^N}\frac{F(v_n)}{v_n^2}w_n^2dx.$$
    Then,
    \begin{equation}
    \begin{array}{lll}
    \frac{1}{2} + o_n(1) & = & \displaystyle\int_{\mathbb{R}^N}\displaystyle\frac{F(v_n)}{v_n^2}w_n^2dx\\
    & \geq & \displaystyle\int_{B_R(y_n)}\displaystyle\frac{F(v_n)}{v_n^2}w_n^2dx\\
    & = & \displaystyle\int_{B_R(0)}\displaystyle\frac{F(v_n(x+y_n))}{v_n(x+y_n)^2}\bar{w}_n^2dx.
    \end{array} \label{eq8}
    \end{equation}
But $v_n(x + y_n) = \|v_n\|^2_\epsilon\bar{w}_n(x) \rightarrow +\infty$ almost everywhere in $\Sigma$. Then, Fatou's Lemma and  $(f_4)$ imply that
    $$\liminf_{n\to\infty}\int_{B_R(0)}\frac{F(v_n(x+ y_n))}{v_n^2(x+ y_n)}\bar{w}^2_n(x)dx = +\infty,$$ which contradicts (\ref{eq8}), and the proof is complete.
\quad\hfill $\Box$

The following result establishes  the existence of a ground-state solution to the corresponding problem to (\ref{P2}) for the case of a constant potential $V$.  The proof can be carried out following the same arguments employed by Rabinowitz in \cite[Theorem 4.23]{Rabinowitz}.

\begin{lemma}
Suppose that $f$ satisfies $(f_1) - (f_5)$. Then, there exists a ground-state solution to the following problem
\begin{equation}
\left\{
\begin{array}{cl}
\Delta^2 w + \alpha w  =  f(w)& \mbox{in $\mathbb{R}^N$}\\
w \in H^2(\mathbb{R}^N),
\end{array}
\right.
\label{P4}
\end{equation}
at the level
$$
c_\alpha = \inf_{\gamma\in\Gamma_\alpha}\sup_{0\leq t\leq 1}I_\alpha(\gamma(t)),$$
where $I_\alpha$ is the energy functional associated to (\ref{P4}) and $$\Gamma_\alpha = \{\gamma\in C^0([0,1],H^2(\mathbb{R}^N)); \,  \gamma(0)=0 \ \ \mbox{and} \ \ I_\alpha(\gamma(1))<0\}.
$$
\label{lemma2.3}
\end{lemma}
The following result gives us an estimate to the   energy level $c_\epsilon$, provided $\epsilon$ is sufficiently small.
\begin{proposition}
Let $c_{V_\infty}$ be the  minimax energy level associated to (\ref{P4}) with $\alpha = V_\infty$.  Then there exists $\epsilon_0 > 0$ such that $c_\epsilon < c_{V_\infty},$ for all $\epsilon \in (0,\epsilon_0)$.
\label{prop2.1}
\end{proposition}
\noindent \textbf{Proof.} Let $w$ be a solution of (\ref{P4}) such that $I_{V_\infty}(w) = c_{V_\infty}$. Fix a function $\chi_R\in
C^1(\mathbb{R}^N,\mathbb{R})$ such that $0\leq\chi_R\leq 1$, $\chi_R
= 1$ in $B_R(0)$, $\chi_R = 0$ in $\mathbb{R}^N \backslash B_{R+2}(0)$ and $|\nabla\chi_R| \leq 1$ in $B_{R+2}(0)\backslash B_R(0)$. Define $v_R(x) = \chi_R(x) w(x)$.
By $(f_2)$ and $(f_3)$,  for all $\eta>0$  there exists $A_\eta > 0$ such that $|f(s)|\leq \eta|s| +
A_\eta|s|^p$. Since  $I_\epsilon'(\varphi_\epsilon(v_R)v_R)\varphi_\epsilon(v_R)v_R = 0$, it follows that
\begin{eqnarray}
\varphi_\epsilon(v_R)^2\int_{\mathbb{R}^N}\left(|\Delta v_R|^2 +
V_\epsilon v_R^2\right)dx &\leq& \eta \varphi_\epsilon(v_R)^2\int_{\mathbb{R}^N}v_R^2dx  \\
&+& A_\eta\varphi_\epsilon(v_R)^{p+1}\int_{\mathbb{R}^N}|v_R|^{p+1}dx.\nonumber
\label{eq9}
\end{eqnarray}
For $\eta = \frac{V_0}{2}$, the previous inequality and $(V_2)$ imply that
\begin{equation}
\varphi_\epsilon(v_R)^2\int_{\mathbb{R}^N}\left(|\Delta v_R|^2 +
\frac{V_0}{2}v_R^2\right)dx \leq C\varphi_\epsilon(v_R)^{p+1}\int_{\mathbb{R}^N}|v_R|^{p+1}dx.
\label{eq10}
\end{equation}
Note that there exists $R_1 > 0$ such that for all $R > R_1$,
\begin{equation}
\int_{B_{R+2}(0)}v_R^2dx \geq \frac{1}{2}\int_{\mathbb{R}^N}w^2dx
\label{eq11}
\end{equation}
and
\begin{equation}
\frac{1}{2}\int_{\mathbb{R}^N}\left(|\Delta w|^2 +
\frac{V_0}{2}w^2\right)dx \leq \int_{B_R}\left(|\Delta w|^2 + \frac{V_0}{2}w^2\right)dx.
\label{eq12}
\end{equation}
From (\ref{eq10}), (\ref{eq12}) and the definition of $v_R$, it follows that
$$\frac{1}{2}\varphi_\epsilon(v_R)^2\int_{\mathbb{R}^N}\left(|\Delta w|^2+\frac{V_0}{2}
w^2\right)dx \leq
C\varphi_\epsilon(v_R)^{p+1}\int_{\mathbb{R}^N}|w|^{p+1}dx,$$ which implies that
\begin{equation}
\varphi_\epsilon(v_R) \geq \left[\frac{\frac{1}{2}\int_{\mathbb{R}^N}\left(|\Delta
w|^2+\frac{V_0}{2}w^2\right)dx}{C\int_{\mathbb{R}^N}|w|^{p+1}dx} \right]^\frac{1}{p-1} =: K > 0, \ \ \forall R>R_1.
\label{eq13}
\end{equation}		
Let $\gamma_R:= \max_{t \geq 0}I_{V_\infty}(tv_R)$ and note that
$$\gamma_R \geq I_{V_\infty}(\varphi_\epsilon(v_R)v_R) = I_\epsilon(\varphi_\epsilon(v_R)v_R) +
\frac{1}{2}\int_{B_{R+2}(0)}\left(V_\infty -
V_\epsilon(x)\right)\varphi_\epsilon(v_R)^2v_R^2dx,$$ which implies that
\begin{equation}
\gamma_R \geq c_\epsilon + \frac{1}{2}\int_{B_{R+2}(0)}\left(V_\infty -
V_\epsilon(x)\right)\varphi_\epsilon(v_R)^2v_R^2dx.
\label{eq14}\end{equation}
We now verify that $\gamma_R = c_{V_\infty} + \psi(R)$, where $\psi(R)\rightarrow 0$ as $R\rightarrow\infty$. In fact, we observe that
 $\gamma_R = I_{V_\infty}(\varphi_{V_\infty}(v_R)v_R) = c_{V_\infty} + I_{V_\infty}(\varphi_{V_\infty}(v_R)v_R) -
I_{V_\infty}(w)$. Following \cite{Weth}, we can prove that $\varphi_{V_\infty}:H^2(\mathbb{R}^N)\rightarrow
\mathbb{R}_+$ is a continuos  function. Since $w$ is a solution of (\ref{P4}) with $\alpha = V_\infty$, it follows that $\varphi_{V_\infty}(v_R)v_R\rightarrow \varphi_{V_\infty}(w)w = w$ as $R\rightarrow \infty$. Therefore
$$\psi(R) = I_{V_\infty}(\varphi_{V_\infty}(v_R)) - I_{V_\infty}(w)\rightarrow 0, \quad \mbox{as $R\rightarrow\infty$.}$$
Take  $R_2>0$ sufficiently large such that
\begin{equation}\psi(R) < \frac{1}{8}\left(V_\infty -
V(0)\right)K^2\int_{\mathbb{R}^N}w^2dx, \quad \mbox{for all $R > R_2$.}\label{eq15}
\end{equation}
As $V_\infty - V(0) > 0$, by the continuity of $V_\infty - V(\cdot)$ in $0$ there exists $\delta > 0$ such that for all $\epsilon < \frac{\delta}{R+2}$, $V_\infty - V(\epsilon x) > \frac{1}{2}(V_\infty - V(0))$ for all $x\in B_{R+2}(0)$. Hence, (\ref{eq14}) implies that
\begin{equation}
\gamma_R \geq c_\epsilon + \frac{1}{4}\int_{B_{R+2}(0)}\left(V_\infty
- V(0)\right)\varphi_\epsilon(v_R)^2v_R^2dx,\label{eq16}
\end{equation}
provided   $0 <\epsilon < {\delta}/{(R+2)}$. Consequently, if $R>R_0:=\max\{R_1,R_2\}$, it follows from (\ref{eq15}), (\ref{eq16}), (\ref{eq13}) and (\ref{eq11}), that
$$c_{V_\infty} + \frac{1}{8}(V_\infty - V(0))K^2\int_{\mathbb{R}^N} w^2dx
> c_\epsilon + \frac{1}{8}(V_\infty - V(0))K^2\int_{\mathbb{R}^N}
w^2dx,$$ which implies that $c_\epsilon < c_{V_\infty},$ provided  $0 < \epsilon <
{\delta}/{(R_0+2)}:=\epsilon_0$.
\quad\hfill $\Box$

\begin{remark}
We observe that if the functional $I_\epsilon$ satisfies the $(PS)_c$ condition for all $c < c_{V_\infty}$, then  the proof of Theorem \ref{theorem2.1} would be complete.
In fact, combing this condition with the Mountain Pass Theorem, there exists a weak nontrivial solution $v_\epsilon$ of (\ref{P3}) for every $\epsilon >0$ sufficiently small.
\label{remark2.1}
\end{remark}

The remainder of this section will be devoted to the proof of the following result.

\begin{proposition}
The functional $I_\epsilon$ satisfies the $(PS)_c$ condition for every $c < c_{V_\infty}$.
\label{proposition2.2}
\end{proposition}

The proof is carried out by a  sequence of lemmas.

\begin{lemma}
Let $\epsilon > 0$ and $(u_n)$ be a $(PS)_c$ sequence for $I_\epsilon$ in $E_\epsilon$, such that $u_n \rightharpoonup u$ in $E_\epsilon$.
The sequence  $v_n:= u_n - u$ is a $(PS)_d$ sequence to $I_\epsilon$, where $d = c - I_\epsilon(u)$.
\label{lemma2.4}
\end{lemma}
\noindent \textbf{Proof.}
We first show that $I_\epsilon(v_n) \rightarrow c - I_\epsilon(u)$, as $n\rightarrow \infty$. In fact, by the weak convergence and Brezis-Lieb Lemma (see \cite{BrezisLieb}), it follows that
\begin{eqnarray*}
I_\epsilon(v_n) - I_\epsilon(u_n) + I_\epsilon(u) & = & \frac{1}{2}\int_{\mathbb{R}^N}\left(|\Delta u_n - \Delta u|^2 - |\Delta u_n|^2 + |\Delta u|^2\right.\\
&  + & \left.  V_\epsilon(x)\left(|u_n - u|^2 - u_n^2 + u^2\right)\right)dx\\
&  - & \int_{\mathbb{R}^N}\left(F(u_n - u) - F(u_n) + F(u)\right)dx\\
& = & \langle u,u \rangle_\epsilon - \langle u_n,u \rangle_\epsilon + o_n(1)\\
& = & o_n(1),
\end{eqnarray*}
and $I_\epsilon(v_n) \rightarrow c - I_\epsilon(u)$, as $n\rightarrow \infty$ as desired.   In order to prove that  $\|I'_\epsilon(v_n)\|_{E_\epsilon^*} = o_n(1)$, from the  weak convergence and the
Brezis-Lieb Lemma, it follows  that
\begin{eqnarray*}
I'_\epsilon(v_n)\varphi - I'_\epsilon(u_n)\varphi & = & \int_{\mathbb{R}^N}\left((\Delta u_n - \Delta u)\Delta\varphi - \Delta u_n\Delta\varphi \right) dx\\
& + & \int_{\mathbb{R}^N}V_\epsilon(x)\left((u_n - u)\varphi - u_n\varphi \right)dx\\
& -  &  \int_{\mathbb{R}^N}\left(f(u_n - u)\varphi - f(u_n)\varphi\right)dx\\
& = & - \langle u,\varphi \rangle_\epsilon +\int_{\mathbb{R}^N}f(u)\varphi dx + o_n(1)\\
& = & I'_\epsilon(u)\varphi + o_n(1) = o_n(1),
\end{eqnarray*}
for every $\varphi \in E_\epsilon$. Therefore, $(v_n)$ is a $(PS)_{c-I_\epsilon(u)}$ sequence.
\quad\hfill $\Box$

\begin{lemma}
Let $\epsilon > 0$ and $(v_n)$ be a $(PS)_d$ sequence to $I_\epsilon$ in $E_\epsilon$. If $v_n \rightharpoonup 0$ in $E_\epsilon$ and $v_n \nrightarrow 0$ in $E_\epsilon$, then
$$c_{V_\infty} \leq d.$$
\label{lemma2.5}
\end{lemma}
\noindent \textbf{Proof.}
Let $s_n > 0$ be such that $s_nv_n \in \mathcal{N}_{V_\infty}$.  We claim that
\begin{equation}
\limsup_{n\to\infty} s_n \leq 1.
\label{eq17}
\end{equation}
In fact, suppose by contradiction that there exist a subsequence $(s_n)$ and $\delta>0$ such that
\begin{equation}
s_n \geq 1 + \delta, \quad \forall\,  n\in \mathbb{N}.
\label{eq18}
\end{equation}
Using the facts that $I_\epsilon'(v_n)v_n = o_n(1)$ and $I_{V_\infty}'(s_n v_n)s_n v_n = 0$ for all $n\in \mathbb{N}$, it follows that
$$\int_{\mathbb{R}^N} \left(\frac{f(s_nv_n)v_n^2}{s_n v_n} - \frac{f(v_n)v_n^2}{v_n}\right)dx = \int_{\mathbb{R}^N}\left(V_\infty - V_\epsilon(x)\right)v_n^2 dx + o_n(1).$$
From $(V_2)$ it follows that for a given $\eta>0$, there exists $R > 0$ such that $V(\epsilon x) \geq V_\infty - \eta$ for all $x \in \mathbb{R}^N$ such that $|x| \geq R\epsilon^{-1}$. Hence,
\begin{eqnarray*}
\int_{\mathbb{R}^N} \left(\frac{f(s_nv_n)v_n^2}{s_n v_n} - \frac{f(v_n)v_n^2}{v_n}\right)dx &\leq &\int_{B_{R\epsilon^{-1}(0)}}\left(V_\infty - V_\epsilon(x)\right)v_n^2 dx \\
& + &   \eta \int_{B_{R\epsilon^{-1}(0)}} v_n^2dx + o_n(1).
\end{eqnarray*}
By  Lemma \ref{lemma2.2}  and the Sobolev embeddings, it follows that
\begin{equation}
\int_{\mathbb{R}^N} \left(\frac{f(s_nv_n)}{s_n v_n} - \frac{f(v_n)}{v_n}\right) v_n^2 dx \leq \eta C + o_n(1).
\label{eq19}
\end{equation}
We now  claim that there exist $R_1,\beta > 0$ and a sequence $(y_n) \subset \mathbb{R}^N$ such that
\begin{equation}
\liminf_{n \to \infty} \int_{B_{R_1}(y_n)} v_n^2 dx \geq \beta.
\label{eq20}
\end{equation}
In fact, on the contrary,  for all $R_1>0$,
$$\lim_{n \to \infty}\sup_{y\in\mathbb{R}^N}\int_{B_{R_1}(y)}v_n^2dx = 0.$$ By Lion's Lemma (see \cite{Lions}), $v_n\rightarrow 0$ as $n\rightarrow
\infty$ in $L^q(\mathbb{R}^N)$ for all $2 < q < 2_*$. Since by $(f_2)$ and $(f_3)$, for $\nu > 0$, there exists $C_\nu > 0$ such that $|f(s)s| \leq \nu|s|^2 + C_\nu|s|^{p+1}$, for all $s\in \mathbb{R}$, it follows from Sobolev embeddings that
$$0 \leq \|v_n\|_\epsilon^2 = \int_{\mathbb{R}^N}f(v_n)(v_n) dx + o_n(1) \leq \nu \int_{\mathbb{R}^N}|v_n|^2dx + C_\nu\int_{\mathbb{R}^N}|v_n|^{p+1}dx + o_n(1).$$ This implies that  $v_n \rightarrow 0$ in $E_\epsilon$, which contradicts our assumption.

Let $\bar{v}_n(x) = v_n(x+y_n)$ and note that using the same arguments that in Lemma \ref{lemma2.2}, one can prove that $(\bar{v}_n)$ is a bounded sequence in $E_\epsilon$. Hence,  $\bar{v}_n \rightharpoonup \bar{v}$ in $E_\epsilon$  along a subsequence. By (\ref{eq20}), $\bar{v} \neq 0$ in a positive measure subset $\Lambda \subset B_{R_1}(0)$. Using Fatou's Lemma, $(f_5)$, (\ref{eq18}) and (\ref{eq19}) it follows that
$$
0 < \int_{\Lambda} \left(\frac{f((1+\delta)\bar{v})}{(1+\delta)\bar{v}} - \frac{f(\bar{v})}{\bar{v}}\right)\bar{v}^2 dx \leq \eta C,
$$
which is impossible because $\forall \eta >0$ is arbitrary.  This contradiction proves that (\ref{eq17}) holds.  Therefore, we have two cases to consider:
\begin{description}
\item [{\it i)}] $\displaystyle \lim_{n\to +\infty} s_n = s < 1$;
\item [{\it ii)}] $\displaystyle \lim_{n\to +\infty}s_n = 1$.
\end{description}

If {\it i)} occurs, then there exists a subsequence $(s_n)$ such that $s_n \rightarrow s < 1$. We can also consider that $s_n < 1$ for all $n\in \mathbb{N}$. From Remark \ref{remark1}, it follows that
\begin{eqnarray*}
c_{V_\infty} & \leq & I_{V_\infty}(s_nv_n)\\
& = & I_{V_\infty}(s_nv_n) - \frac{1}{2}I'_{V_\infty}(s_nv_n)s_nv_n\\
& = & \int_{\mathbb{R}^N}\frac{1}{2}\left(f(s_nv_n)s_nv_n - 2F(s_nv_n)\right)dx\\
& \leq & \int_{\mathbb{R}^N}\frac{1}{2}\left(f(v_n)v_n - 2F(v_n)\right)dx\\
& = & I_\epsilon(v_n) - \frac{1}{2}I'_\epsilon(v_n)v_n + o_n(1)\\
& = & d + o_n(1).
\end{eqnarray*}
Taking $n \to \infty$, we obtain that $c_{V_\infty} \leq d$ as required.
which is the desired conclusion.

Suppose that {\it ii)} holds. In this case,
    $$d + o_n(1) = I_\epsilon(v_n) = I_{V_\infty}(s_nv_n) + I_\epsilon(v_n) - I_{V_\infty}(s_nv_n),$$ which implies that
    \begin{equation}
d + o_n(1) = I_\epsilon(v_n) \geq c_{V_\infty} + I_\epsilon(v_n) - I_{V_\infty}(s_nv_n).
\label{eq21}
\end{equation}
Therefore, it remains to prove that $I_\epsilon(v_n) - I_{V_\infty}(s_nv_n) = o_n(1)$.
Note that
\begin{eqnarray}
I_\epsilon(v_n) - I_{V_\infty}(s_nv_n) & = & \int_{\mathbb{R}^N}\frac{(1 - s_n^2)}{2}|\Delta v_n|^2dx + \frac{1}{2}\int_{\mathbb{R}^N} V_\epsilon(x)v_n^2dx  \\
& -&   \frac{s_n^2}{2}\int_{\mathbb{R}^N}V_\infty v_n^2dx + \int_{\mathbb{R}^N}\left(F(s_n v_n) - F(v_n)\right)dx.\nonumber
\label{eq22}
\end{eqnarray}
Since  $(v_n)$is bounded  in $E_\epsilon$,
$$\int_{\mathbb{R}^N}\frac{(1 - s_n^2)}{2}|\Delta v_n|^2dx = o_n(1).$$
For any $R > 0$,  the Sobolev embeddings and the continuity of $V$ imply
$$\frac{1}{2}\int_{\overline{B_R(0)}} V_\epsilon(x)v_n^2dx = o_n(1)$$
and
$$\frac{s_n^2}{2}\int_{\overline{B_R(0)}} V_\infty v_n^2dx = o_n(1).$$
Hence,
\begin{eqnarray*}
I_\epsilon(v_n) - I_{V_\infty}(s_nv_n) & \geq & o_n(1) +  \frac{1}{2}\int_{B_R(0)^c} V_\epsilon(x)v_n^2dx - \frac{s_n^2}{2}\int_{B_R(0)^c}V_\infty v_n^2dx\\
& & + \int_{\mathbb{R}^N}\left(F(s_n v_n) - F(v_n)\right)dx.
\end{eqnarray*}
By $(V_2)$, given $\eta > 0$  there exists $R>0$ sufficient large  such that
\begin{eqnarray*}
I_\epsilon(v_n) - I_{V_\infty}(s_nv_n) & \geq & o_n(1) +  \frac{1}{2}\int_{B_R(0)^c} (V_\infty - \eta) v_n^2dx - \frac{s_n^2}{2}\int_{B_R(0)^c}V_\infty v_n^2dx\\
& & + \int_{\mathbb{R}^N}\left(F(s_n v_n) - F(v_n)\right)dx,
\end{eqnarray*}
which implies that
\begin{eqnarray*}
I_\epsilon(v_n) - I_{V_\infty}(s_nv_n) & \geq & o_n(1) +  \frac{(1 - s_n^2)}{2}\int_{B_R(0)^c} V_\infty v_n^2dx - \frac{\eta}{2}\int_{B_R(0)^c}v_n^2dx\\
& & + \int_{\mathbb{R}^N}\left(F(s_n v_n) - F(v_n)\right)dx.
\end{eqnarray*}
Using that  $(v_n)$  is bounded  and  Sobolev embeddings, yields
\begin{equation}
I_\epsilon(v_n) - I_{V_\infty}(s_nv_n) \geq o_n(1) - C \eta + \int_{\mathbb{R}^N}\left(F(s_n v_n) - F(v_n)\right)dx.
\label{eq23}
\end{equation}
By (\ref{eq23}) and (\ref{eq21}), we have
$$d + o_n(1) \geq c_{V_\infty} - C\eta + o_n(1) + \int_{\mathbb{R}^N}\left(F(s_n v_n) - F(v_n)\right)dx.$$
By the mean value theorem,  $\int_{\mathbb{R}^N}\left(F(s_n v_n) - F(v_n)\right)dx = o_n(1)$. Thus,
$$d + o_n(1) \geq c_{V_\infty} - C\eta + o_n(1)$$ and the result follows after passing to the limit  $n \to \infty$. \quad\hfill $\Box$

As a consequence of the above lemma, we have:
\begin{corollary}
If $(v_n)$ is a $(PS)_d$ sequence for $I_\epsilon$ such that $v_n \rightharpoonup 0$ and $d < c_{V_\infty}$, then $v_n \rightarrow 0$ in $E_\epsilon$.
\label{corollary2.1}
\end{corollary}

Finally we can proceed with the proof of Propostion \ref{proposition2.2}.

\medskip
\noindent \textbf{Proof of Proposition \ref{proposition2.2}.}  Let $(u_n)$ be a $(PS)_c$ sequence for $I_\epsilon$. By Lemma \ref{lemma2.2}, $(u_n)$ is a bounded sequence in $E_\epsilon$. Then there exists $u\in E_\epsilon$ such that $u_n \rightharpoonup u$ in $E_\epsilon$. If we denote by $v_n = u_n - u$, it follows that $v_n \rightharpoonup 0$ in $E_\epsilon$. By Lemma \ref{lemma2.4}, it follows that $(v_n)$ is a  $(PS)_d$ sequence for $I_\epsilon$, where $d = c - I_\epsilon(u)$.
Since $u$ is a weak solution of (\ref{P3}), then $I_\epsilon(u) \geq c_\epsilon > 0$ and $d \leq c < c_{V_\infty}$. By Corollary \ref{corollary2.1},
$v_n\rightarrow 0$ in $E_\epsilon$ and proof is complete. \quad\hfill $\Box$

Therefore, Remark \ref{remark2.1} implies that there exists a nontrivial weak solution $v_\epsilon$ to (\ref{P3}) for every $\epsilon >0$ sufficiently small, and Theorem \ref{theorem2.1} follows.

\section{Concentration}

\noindent In this section our goal is to prove the concentration phenomenon stated in Theorem \ref{theorem1.1}.   Invoking  Lemma \ref{lemma2.3}, let $w \in H^2(\mathbb{R}^N)$ be a ground state solution  to the problem
\begin{equation}
\left\{
\begin{array}{cl}
\Delta^2 u + V_0u  =  f(u) & \,\ \mbox{em $\mathbb{R}^N$}\\
u \in H^2(\mathbb{R}^N),
\end{array}\right. \label{P5}
\end{equation}
We begin by showing the following limit:

\begin{lemma}
$$\lim_{\epsilon \to 0} c_\epsilon = c_{V_0}.$$
\label{lemma3.1}
\end{lemma}

\noindent \textbf{Proof.}   Let $\psi\in C^\infty_0(\mathbb{R}^N)$ be such that $0\leq
\psi \leq 1$, $\psi \equiv 0$ in $\mathbb{R}^N\slash B_2(0)$, $\psi \equiv 1$ in
$B_1(0)$, $|\nabla \psi|\leq C$ and $|\Delta \psi|\leq C$ in $\mathbb{R}^N$. Let us define
$$w_\epsilon(x) = \psi(\epsilon x)w(x).$$
Note that $w_\epsilon \rightarrow w$ in $H^2(\mathbb{R}^N)$ and $I_{V_0}(w_\epsilon) \rightarrow I_{V_0}(w)$ as $\epsilon\rightarrow 0$ where $I_{V_0}$ is the energy functional associated to (\ref{P5}).
Let $\varphi_\epsilon(w_\epsilon)$ be such that $\varphi_\epsilon(w_\epsilon)w_\epsilon \in \mathcal{N}_\epsilon.$ Suppose that $\varphi_\epsilon(w_\epsilon) \rightarrow 1$ as $\epsilon \rightarrow 0$. Note that
\begin{eqnarray*}
c_\epsilon & \leq & I_\epsilon(\varphi_\epsilon(w_\epsilon)w_\epsilon)\\
& = & I_{V_0}(\varphi_\epsilon(w_\epsilon)w_\epsilon) +
\frac{1}{2}\int_{\mathbb{R}^N}\varphi_\epsilon(w_\epsilon)^2\left(V_\epsilon(\epsilon
x) - V_0\right)w_\epsilon^2 dx.
\end{eqnarray*}
Using the Lebesgue Dominated Theorem, it follows that
$$\limsup_{\epsilon \to 0} c_\epsilon = I_{V_0}(w) = c_{V_0}.$$
On the other hand, since $I_{V_0}(v) \leq I_\epsilon(v)$ for all $v\in H^2(\mathbb{R}^N)$, it follows that $c_{V_0} \leq c_\epsilon$. Then
$$\lim_{\epsilon \to 0} c_\epsilon = c_{V_0}.$$
It remains to prove that $\varphi_\epsilon(w_\epsilon) \rightarrow 1$ as $\epsilon\rightarrow 0$.
Since $I_\epsilon'(\varphi_\epsilon(w_\epsilon)w_\epsilon)w_\epsilon = 0$, it follows that
$$\varphi_\epsilon(w_\epsilon)\int_{\mathbb{R}^N}\left(|\Delta w_\epsilon|^2 + V(\epsilon x)w_\epsilon^2\right)dx = \int_{\mathbb{R}^N}f(\varphi_\epsilon(w_\epsilon)w_\epsilon)w_\epsilon dx. $$
We claim that $(\varphi_\epsilon(w_\epsilon))$ is bounded. In fact, on the contary,
 there exists $\epsilon_n \to 0$ such that $\varphi_{\epsilon_n}(w_{\epsilon_n}) \rightarrow +\infty$. Let $\Sigma \subset \mathbb{R}^N$ be such that $|\Sigma| > 0$ and $w(x) \neq 0$ for all $x\in \Sigma$. Hence, calling Remark \ref{remark1}, it holds for all $n \in \mathbb{N}$ that
\begin{eqnarray*}
\|w_{\epsilon_n}\|_{\epsilon_n}^2 & = & \int_{\mathbb{R}^N}\frac{f(\varphi_{\epsilon_n}(w_{\epsilon_n})w_{\epsilon_n})\varphi_{\epsilon_n}(w_{\epsilon_n})w_{\epsilon_n}}{\varphi_{\epsilon_n}(w_{\epsilon_n})^2}dx\\
& \geq &  \int_{\Sigma}\frac{2F(\varphi_{\epsilon_n}(w_{\epsilon_n})w_{\epsilon_n})}{\varphi_{\epsilon_n}(w_{\epsilon_n})^2}dx\\
& = &  \int_{\Sigma \backslash w_{\epsilon_n}^{-1}(0)^c}\frac{2F(\varphi_{\epsilon_n}(w_{\epsilon_n})w_{\epsilon_n})}{(\varphi_{\epsilon_n}(w_{\epsilon_n})w_{\epsilon_n})^2}w_{\epsilon_n}^2dx.
\end{eqnarray*}
On the other hand, by $(f_4)$ and Fatou's Lemma it follows that
\begin{eqnarray*}
\liminf_{n \rightarrow \infty} \int_{\Sigma \backslash w_{\epsilon_n}^{-1}(0)^c}\frac{2F(\varphi_{\epsilon_n}(w_{\epsilon_n})w_{\epsilon_n})}{(\varphi_{\epsilon_n}(w_{\epsilon_n})w_{\epsilon_n})^2}w_{\epsilon_n}^2dx = +\infty,
\end{eqnarray*}
which implies that
$$\|w_{\epsilon_n}\|_{\epsilon_n}^2 \rightarrow +\infty, \quad \mbox{as $n\to\infty$},$$ which contradicts the fact that $w_{\epsilon_n} \rightarrow w$ as $n \rightarrow \infty$.

We can now verify that $\varphi_\epsilon(w_\epsilon) \nrightarrow 0$ as $\epsilon\rightarrow 0$. In fact, on the contrary there exists $\epsilon_n \rightarrow 0$ such that  $\varphi_{\epsilon_n}(w_{\epsilon_n}) \rightarrow 0$ as $n\rightarrow\infty$. By $(f_2)-(f_3)$ one can prove that
  \begin{equation}\label{limite1}
  \lim_{n\to \infty}\int_{\mathbb{R}^N}\frac{f(\varphi_{\epsilon_n}(w_{\epsilon_n})w_{\epsilon_n})w_{\epsilon_n}^2}{\varphi_{\epsilon_n}(w_{\epsilon_n})w_{\epsilon_n}} dx = 0.
  \end{equation}
 On the other hand,
 \begin{equation}\label{limite2}
 \|w_{\epsilon_n}\|_{\epsilon_n}^2 = \int_{\mathbb{R}^N}\frac{f(\varphi_{\epsilon_n}(w_{\epsilon_n})w_{\epsilon_n})w_{\epsilon_n}^2}{\varphi_{\epsilon_n}(w_{\epsilon_n})w_{\epsilon_n}} dx
 \end{equation}
Hence by (\ref{limite1}) and (\ref{limite2}), one can see that $\|w_{\epsilon_n}\|_{\epsilon_n} \to 0$, which contradicts the fact that $w_{\epsilon_n }\to w$ and $I_{V_0}(w) = c_{V_0} >0$.
Then there exist $\alpha,\beta >0$ such that
$$\alpha \leq \varphi_{\epsilon}(w_{\epsilon}) \leq \beta.$$
Using that $w_{epsilon_n} \to w$ in $H^2(\mathbb{R}^N)$ and $w$ is a solution of (\ref{P5}), it follows by $(f_5)$ that
 $\varphi_\epsilon(w_\epsilon) \rightarrow 1$.
\quad\hfill $\Box$

In the following, we consider  a sequence $(\epsilon_n)$, with $\epsilon_n \rightarrow 0$ as $n\to \infty$,  and let $u_{\epsilon_n}$  be a solution of  (\ref{P1}) given by Theorem \ref{theorem2.1}. Let $v_n(x) := v_{\epsilon_n}(x) = u_{\epsilon_n}(\epsilon_n x)$ . Similar arguments employed in  proof of Lemma \ref{lemma2.2} show that $(v_n)$ is a bounded sequence in $H^2(\mathbb{R}^N)$.

\begin{lemma}
There exists $(y_n)\subset \mathbb{R}^N$ and $R,\beta >0$ such that
$$\liminf_{n \to \infty} \int_{B_R(y_n)} v_n^2 dx \geq \beta > 0.$$ \label{lemma5}
\end{lemma}
\noindent \textbf{Proof.}
Suppose the assertion of the lemma is false. Then by Lemma I.1 of \cite{Lions} (with $q=2$ and $p=\frac{2N}{N-2}$ ), $v_n\rightarrow 0$ in $L^r(\mathbb{R}^N)$ where $2 < r < 2_*$. Hence by the Lebesgue Dominated Convergence Theorem, we get
$$\int_{\mathbb{R}^N} f(v_n) v_n dx = o_n(1) \,\ \mbox{and} \,\ \int_{\mathbb{R}^N} F(v_n)dx = o_n(1).$$
Then $c_{\epsilon_n}\rightarrow 0$ as $n\rightarrow \infty$, which contradicts Lemma \ref{lemma3.1} and this contradiction proves the lemma.
\quad\hfill $\Box$

Define the  function $w_n(x) = v_n(x+y_n) = u_n(\epsilon_n x + \epsilon_n y_n)$. Note that $w_n$ satisfies
\begin{equation}
\left\{
\begin{array}{lll}
\Delta^2w_n + V(\epsilon_n x + \epsilon_n y_n)w_n & = &
f(w_n) \,\ \mbox{in $\mathbb{R}^N$}\\
w_n \in H^2(\mathbb{R}^N),
\end{array} \right. \label{P6}
\end{equation}
and
\begin{equation}
\liminf_{n\rightarrow \infty}\int_{B_R(0)}w_n^2 dx \geq
\beta. \label{eq24}
\end{equation}

\begin{lemma}
The sequence $(\epsilon_ny_n)$ is bounded  in $\mathbb{R}^N$.
\end{lemma}
\noindent \textbf{Proof.}
Suppose by contradiction that there exists a subsequence $(\epsilon_n y_n)$ such that $\epsilon_n y_n\rightarrow \infty$ as $n\to \infty$. Since $(w_n)$ is a bounded sequence, there exists $w_0 \in H^2(\mathbb{R}^N)$ such that $w_n\rightharpoonup w_0$ in $H^2(\mathbb{R}^N)$ and $w_n\rightarrow w_0$ in $L^q_{loc}(\mathbb{R}^N)$ where $2 \leq q < 2_*$ as $n\rightarrow \infty$. Note that by (\ref{eq24}), $w_0 \neq 0$. From $(V_2)$, one can prove that $w_0$ satisfies (\ref{P4}) with $\alpha = V_\infty$.

Using $(V_2)$, Lemma \ref{lemma3.1} and Fatou's Lemma, we get
\begin{eqnarray*}
c_{V_0} & < & c_{V_\infty}\\
& \leq & I_{V_\infty}(w_0) - \frac{1}{2}I'_{V_\infty}(w_0)w_0\\
& = & \int_{\mathbb{R}^N}\left(\frac{1}{2}f(w_0)w_0 - F(w_0)\right)dx\\
& \leq & \liminf_{n\to\infty}\int_{\mathbb{R}^N}\left(\frac{1}{2}f(w_n)w_n - F(w_n)\right)dx\\
& = & \liminf_{n\to\infty} c_{\epsilon_n} = c_{V_0},
\end{eqnarray*}
which give us a contradiction.

\quad\hfill $\Box$

Note that by the last result, we can assume that there exists $x_0 \in \mathbb{R}^N$ such that $\epsilon_n y_n \rightarrow x_0$ as $n\rightarrow \infty$. We can suppose also that $w_n \rightharpoonup w_0$ in $H^2(\mathbb{R}^N)$ where $w_0 \neq 0$.

\begin{lemma}
The point $x_0$ is a global minimum to $V$.
\label{lemma3.4}
\end{lemma}

\noindent \textbf{Proof.}
By (\ref{P6}) and elliptic regularity theory, one can prove that in fact $w_n\rightarrow w_0$ in $C^4_{loc}(\mathbb{R}^N)$ as $n\rightarrow \infty$. Then for each $x\in\mathbb{R}^N$, $w_0$ satisfies the following equation
$$
\Delta^2 w_0(x) + V(x_0)w_0(x) = f(w_0(x)).
$$
Hence,
$$\int_{\mathbb{R}^N}\left(|\Delta w_0|^2 + V_0w_0^2\right)dx \leq
\int_{\mathbb{R}^N}\left(|\Delta w_0|^2 + V(x_0)w_0^2\right)dx =
\int_{\mathbb{R}^N}f(w_0)w_0dx$$ and there exists $0<\tau\leq 1$
such that $\tau w_0\in\mathcal{N}_{V_0}$, where $\mathcal{N}_{V_0}$ denotes the Nehari manifold associated to (\ref{P5}). Fatou's Lemma and Remark \ref{remark1} imply that
\begin{eqnarray*}
c_{V_0} & = & \lim_{n\rightarrow\infty} c_{\epsilon_n}\\
& = &
\liminf_{n\rightarrow\infty}\int_{\mathbb{R}^N}\left(\frac{1}{2}f(w_n)w_n
dx - F(w_n)\right)dx\\
& \geq & \int_{\mathbb{R}^N}\left(\frac{1}{2}f(w_0)w_0
dx - F(w_0)\right)dx\\
& \geq & \int_{\mathbb{R}^N}\left(\frac{1}{2}f(\tau w_0)\tau w_0
dx - F(\tau w_0)\right)dx\\
& = & I_{V_0}(\tau w_0) \geq c_{V_0}
\end{eqnarray*}
and this implies that $\tau = 1$.
Therefore $w_0\in\mathcal{N}_{V_0}$ and
$$\int_{\mathbb{R}^N}\left(|\Delta w_0|^2 + V(x_0)w_0^2\right)dx =
\int_{\mathbb{R}^N}f(w_0)w_0 dx =
\int_{\mathbb{R}^N}\left(|\Delta w_0|^2 + V_0w_0^2\right)dx,$$ which implies that $V(x_0) = V_0$.
\quad\hfill $\Box$

\begin{lemma}
$w_n \rightarrow w_0$ in $H^2(\mathbb{R}^N)$ as $n \to \infty$.
\label{lemma3.5}
\end{lemma}
\noindent \textbf{Proof.}
By Lemma \ref{lemma3.1}, we have
$$
\lim_{n\to\infty} I_{\epsilon_n}(v_n) = \lim_{n\to\infty} c_{\epsilon_n} = c_{V_0}.
$$
Given $v\in H^2(\mathbb{R}^N)\backslash\{0\}$, from $(f_5)$,  there exists  $\varphi_{V_0}(v)>0$ such that $\varphi_{V_0}(v)v \in \mathcal{N}_{V_0}$. Set $\tilde{w}_n = \varphi_{V_0}(w_n)w_n$. Hence,
\begin{eqnarray*}
c_{V_0} & \leq & \frac{1}{2}\int_{\mathbb{R}^N}\left(|\Delta\tilde{w}_n|^2 + V_0\tilde{w}_n^2\right) dx - \int_{\mathbb{R}^N}F(\tilde{w}_n)dx\\
& \leq & \frac{1}{2}\int_{\mathbb{R}^N}\left(|\Delta\tilde{w}_n|^2 + V(\epsilon_n x + \epsilon_n y_n)\tilde{w}_n^2\right) dx - \int_{\mathbb{R}^N}F(\tilde{w}_n)dx\\
& = & I_{\epsilon_n}(\varphi_{V_0}(w_n)v_n)
\leq I_{\epsilon_n}(v_n)
 = c_{\epsilon_n} = c_{V_0} + o_n(1),
\end{eqnarray*}
which implies that $I_{V_0}(\tilde{w}_n)\rightarrow c_{V_0}$ as $n\rightarrow \infty$.

We now prove that $\varphi_{V_0}(w_n) \rightarrow \varphi_0 > 0$ along a subsequence. We  first observe  that there exists $M>0$ such that $|\varphi_{V_0}(w_n)| \leq M$, $\forall n\in\mathbb{N}$. In fact, since $w_n \nrightarrow 0$ there exists $\delta>0$ such that $\|w_n\|_{H^2(\mathbb{R}^N)} > \delta$ along a subsequence. On the other hand,  since $I_{V_0}(\tilde{w}_n) \rightarrow c_{V_0}$ and $I_{V_0}'(\tilde{w}_n)\tilde{w}_n = 0$ for all $n\in \mathbb{N}$, it is easy to see that $(\tilde{w}_n)$ is a bounded sequence in $H^2(\mathbb{R}^N)$. Then
$$|\varphi_{V_0}(w_n)|\delta < \|\varphi_{V_0}(w_n)w_n\|_{H^2(\mathbb{R}^N)} \leq K$$
which implies that
$$|\varphi_{V_0}(w_n)| \leq \frac{K}{\delta} = M, \ \ \forall n\in\mathbb{N}.$$
Hence,  $\varphi_{V_0}(w_n)\rightarrow \varphi_0 \geq 0$.  We now observe that $\varphi_0 > 0$, otherwise
$$\|\tilde{w}_n\|_{H^2(\mathbb{R}^N)} = |\varphi_{V_0}(w_n)|\|w_n\|_{H^2(\mathbb{R}^N)} \rightarrow 0$$ as $n\rightarrow \infty$, which is impossible.  Therefore $\tilde{w}_n = \varphi_0(w_n)w_n \rightharpoonup\varphi_0 w \neq 0$ in $H^2(\mathbb{R}^N)$.  Therefore, we conclude the lemma from the next result. \quad\hfill $\Box$

In the proof of the next result we use some arguments of Alves and Figueiredo found in \cite{Alves1}.

\begin{lemma}
Let $(z_n) \subset H^2(\mathbb{R}^N)$ be a sequence such that $I_{V_0}(z_n) \rightarrow c_{V_0}$ as $n\to \infty$ and $z_n \in \mathcal{N}_{V_0}$ for all $n\in\mathbb{N}$. If $z_n \rightharpoonup z \neq 0$ in $H^2(\mathbb{R}^N)$, then $z_n \rightarrow z$ in $H^2(\mathbb{R}^N)$ along a subsequence.
\label{lemma3.6}
\end{lemma}
\noindent \textbf{Proof.}

By the Ekeland Variational Principle, we can assume that $(z_n)$ is a $(PS)_{c_{V_0}}$ sequence for $I_{V_0}$ in $H^2(\mathbb{R}^N)$. Then it is possible to show that $I_{V_0}'(z) = 0$ which implies that $z\in \mathcal{N}_{V_0}$.

Using Remark \ref{remark1} and Fatou's Lemma, it follows that
\begin{eqnarray*}
c_{V_0} & = & \lim_{n\to\infty}\left[ \int_{\mathbb{R}^N}\left(\frac{1}{2}f(z_n)z_n - F(z_n)\right)dx + o_n(1)\right]\\
 & \geq & \int_{\mathbb{R}^N}\left(\frac{1}{2}f(z)z - F(z)\right)dx\\
 & = & I_{V_0}(z)\\
 & \geq & c_{V_0},
\end{eqnarray*}
which implies that
\begin{equation}
I_{V_0}(z) = c_{V_0}. \label{eq25}
\end{equation}

Let $v_n = z_n - z$ and note that by Brezis-Lieb Lemma, $(v_n)$ is $(PS)_d$ sequence for $I_{V_0}$ where $d = c_{V_0} - I_{V_0}(z) = 0$. Note that $v_n \rightharpoonup 0$ in $H^2(\mathbb{R}^N)$ and we claim that in fact $v_n \rightarrow 0$ in $H^2(\mathbb{R}^N)$. On the contrary, if $v_n \nrightarrow 0$ in $H^2(\mathbb{R}^N)$, we can use the same arguments than in Lemma \ref{lemma2.5} to prove that $(v_n)$ is a $(PS)_d$ sequence to $I_{V_0}$ for $d \geq c_{V_0} > 0$. But this contradicts the fact that $(v_n)$ is a $(PS)_0$ sequence and this contradiction proves the lemma.
\quad\hfill $\Box$

\medskip
Combing  Lemma \ref{lemma3.5} with the Sobolev embeddings, it follows that $w_n \rightarrow w$ in $L^{2_*}(\mathbb{R}^N)$ as $n \to \infty$. Therefore,  we obtain
\begin{equation}
\int_{B_R^c(0)}|w_n|^{2_*}dx \rightarrow 0 \quad \mbox{as $R\rightarrow \infty$ uniformly in $n$}.
\label{eq3}
\end{equation}

\begin{lemma}
$w_n(x) \rightarrow 0$ as $|x|\rightarrow \infty$, uniformly in $n$. \label{lemma3.7}
\end{lemma}
\noindent \textbf{Proof.}
By the uniform $L^\infty$ estimates to solutions of subcritical biharmonic equations given by Ramos in \cite{Ramos}, we have
$$\|w_n\|_{L^\infty(\mathbb{R}^N)} \leq C, \quad \forall\, n\in\mathbb{N},$$ where $C$ is independent of $n$.
Given any $x\in\mathbb{R}^N$, the function $w_n \in L^q(B_1(x))$ for all $q\geq 1$. By \cite[Theorem 7.1]{Agmon} it follows that
\begin{eqnarray*}
\|w_n\|_{W^{4,q}(B_1(x))} & \leq & C\left( \|f(w_n)\|_{L^q(B_2(x))} + \|w_n\|_{L^q(B_2(x))}\right) \\
& \leq & C \|w_k\|_{L^q(B_2(x))}\\
& \leq & C \|w_k\|_{L^\infty(\mathbb{R}^N)}^{\frac{q-2_*}{q}} \|w_k\|_{L^{2_*}(B_2(x))}^{2_*}\\
& = & C \|w_k\|_{L^{2_*(B_2(x))}}^{2_*},
\end{eqnarray*}
with $C>0$ being a constant independent of $x$ and $n$.   If $q>N$, we have the continuous imbedding $W^{4,q}(B_1(x))\hookrightarrow C^{3,\alpha}(\overline{B_1(x)})$ for $\alpha \in \left(0,1-\frac{N}{q}\right)$. Then
\begin{eqnarray*}
\|w_k\|_{C^{3,\alpha}(\overline{B_1(x)})} & \leq & \|w_k\|_{W^{4,q}(B_1(x))}
 \leq C \|w_k\|_{L^{2_*}(B_2(x))}^{2_*}.
\end{eqnarray*}
By (\ref{eq3}), it follows that
$|w_n(x)| \rightarrow 0$   as $|x|\rightarrow \infty$ uniformly in $n$.
\quad\hfill $\Box$

In order to prove the concentration  behavior  of solutions, we claim that there exists $\rho>0$ such that $\|u_n\|_{L^\infty(\mathbb{R}^N)} = \|w_n\|_{L^\infty(\mathbb{R}^N)} > \rho$, for all $n\in\mathbb{N}$ along a subsequence. In fact, if $\|w_n\|_{L^\infty(\mathbb{R}^N)} \rightarrow 0$, since for all $\eta > 0$ there exists $A_\eta > 0$ such that $|f(s)s| \leq \eta |s|^2 + A_\eta |s|^{p+1}$, for all $s \in \mathbb{R}$, it follows that
\begin{eqnarray*}
\|w_n\|^2_{H^2(\mathbb{R}^N)} & \leq & C\int_{\mathbb{R}^N}\left(|\Delta w_n|^2 + V(\epsilon_n x + \epsilon_n y_n)w_n\right)dx\\
& = & C\int_{\mathbb{R}^N}f(w_n)w_n dx\\
& \leq & C\left(\eta\|w_n\|^2_{L^2(\mathbb{R}^N)} + A_\eta\|w_n\|^{p+1}_{L^{p+1}(\mathbb{R}^N)}\right).
\end{eqnarray*}
In particular, for $0< \eta < {1}/{2}$, we have
$$
\|w_n\|^2_{H^2(\mathbb{R}^N)} \leq A_\eta\|w_n\|_{L^\infty(\mathbb{R}^N)}^{p+1}\|w_n\|^p_{L^p(\mathbb{R}^N)} \rightarrow 0
$$
Hence, if  $\|w_n\|_{L^\infty(\mathbb{R}^N)} \to 0$, then  $\|w_n\|^2_{H^2(\mathbb{R}^N)} \rightarrow 0$ as $n\to\infty$, which contradicts the fact that $w_n\rightarrow w$ and $w \neq 0$.

Let $x_n$ be the maximum point of $|u_n|$ in $\mathbb{R}^N$.  Then
$$
p_n := \frac{x_n - \epsilon_ny_n}{\epsilon_n}
$$
is the maximum point of $|w_n|$. By Lemma \ref{lemma3.7}, there exists $R_0 > 0$ such that $p_n \in B_{R_0}(0)$ for all sufficiently large $n$. Then, along a subsequence $p_n \rightarrow p_0$ as $n\rightarrow \infty$. Hence
$$
x_n = \epsilon_n p_n + \epsilon_n y_n \rightarrow x_0 \ \ \mbox{as $n\rightarrow \infty$},
$$
which proves Theorem \ref{theorem1.1}.

\noindent \textbf{Acknowledgment.}\  The authors are grateful to  Profs. Claudianor O. Alves and Marco A. S. Souto for valuable discussions.

\end{document}